\documentclass[10pt]{amsart}
\usepackage{amsmath,amsfonts,amssymb}
\usepackage{enumerate}
\newtheorem{lemma}{Lemma}[section]

\newtheorem{theorem}[lemma]{Theorem}

\theoremstyle{definition}

\theoremstyle{remark}

\newcommand{\Mat}{\mathbf{M}_d(\mathbb{K})}

\numberwithin{equation}{section} \numberwithin{table}{section}

\newcommand{\vvv}{{|\!|\!|}}

\begin{document}
\title[Uniqueness of Barabanov norms]{A new sufficient condition for the uniqueness of Barabanov norms}
%\dedicatory{for Stronzo Bestiale of Palermo}
\author{Ian D. Morris}
\address{Dipartimento di Matematica, Universit\'a degli Studi di Roma Tor Vergata, Via della Ricerca Scientifica, Rome 00133, Italy.}
\email{ian.morris.ergodic@gmail.com}

%\abstract{We continue our investigation into the uniqueness of Barabanov norms. We prove a theorem. We provide some examples of the application of the Theorem which exhibit behaviour contrasting with that considered in earlier work, and also give a theoretical application dealing with the prevalence or otherwise of unique Barabanov norms.}

\begin{abstract}
The joint spectral radius of a bounded set of $d \times d$ real or complex matrices is defined to be the maximum exponential rate of growth of products of matrices drawn from that set. Under quite mild conditions such a set of matrices admits an associated vector norm, called a Barabanov norm, which can be used to characterise those sequences of matrices which achieve this maximum rate of exponential growth. In this note we continue an earlier investigation into the problem of determining when the Barabanov norm associated to such a set of matrices is unique. We give a new sufficient condition for this uniqueness, and provide some examples in which our condition applies. We  also give a theoretical application which shows that the property of having a unique Barabanov norm can in some cases be highly sensitive to small perturbations of the set of matrices.\end{abstract}

\maketitle

\section{Introduction}

Given a bounded set $\mathcal{A}$ of $d \times d$ matrices over $\mathbb{R}$ or $\mathbb{C}$, the joint spectral radius of $\mathcal{A}$ is defined to be the quantity
\[\varrho(\mathcal{A})=\lim_{n \to \infty} \sup\left\{\left\|A_{i_n} \cdots A_{i_1}\right\|^{\frac{1}{n}} \colon A_i \in \mathcal{A}\right\},\]
a definition introduced by G.-C. Rota and G. Strang in 1960 (\cite{RS}, subsequently reprinted in \cite{Rotacoll}). This limit always exists and is independent of the norm used (for a proof see e.g. \cite{Jungers}). The joint spectral radius has subsequently been found to arise in a range of mathematical contexts including control and optimisation \cite{Ba,BM,Gu}, wavelet regularity \cite{DL0}, coding theory \cite{MO}, and combinatorics \cite{BCJ}. As such, the properties of the joint spectral radius are the subject of ongoing research investigation (see for example \cite{AAJPR,BN,BM,Chang, GWZ,HMST,Kozb2}). This note is concerned with a theoretical tool associated to the joint spectral radius called the \emph{Barabanov norm}, which we now define.

Let us say that a set $\mathcal{A}$ of $d \times d$ real or complex matrices is \emph{reducible} if its elements simultaneously preserve a linear subspace with dimension strictly between $0$ and $d$. If $\mathcal{A}$ is not reducible then it will be called \emph{irreducible}. An irreducible set of matrices always has nonzero joint spectral radius (see e.g. \cite{Jungers}). In the article \cite{Ba}, N. E. Barabanov showed that to any compact irreducible set $\mathcal{A}$ of $d \times d$ matrices over $\mathbb{R}$ or $\mathbb{C}$, one may associate a norm $\vvv \cdot \vvv$ on $\mathbb{R}^d$ or $\mathbb{C}^d$ such that the Bellman-like equation 
\begin{equation}\label{barab}\varrho(\mathcal{A})\vvv v\vvv=\sup\{\vvv Av\vvv \colon A \in \mathcal{A}\}\end{equation}
is satisfied for every vector $v$. We shall call a norm which satisfies this relation for every $v$ a \emph{Barabanov norm} for $\mathcal{A}$. Given any vector $v$, by iterating the above relation it follows that for each $n \geq 1$,
\[\varrho(\mathcal{A})^n\vvv v\vvv=\sup\{\vvv A_{i_n} \cdots A_{i_1}v\vvv \colon A_i \in \mathcal{A}\},\]
and since $\mathcal{A}$ is assumed to be compact, it follows that we may extract a sequence $(A_{i_j})_{j=1}^\infty$ of elements of $\mathcal{A}$ such that $\vvv A_{i_n}\cdots A_{i_1}v\vvv=\varrho(\mathcal{A})^n \vvv v\vvv$ for every $n \geq 1$. As various researchers have noted (see e.g. \cite{Jungers,Koz3,Marg}), Barabanov norms thus implicitly encode a description of certain sequences of matrices drawn from $\mathcal{A}$ whose partial products grow at the maximum possible exponential rate. The problem of constructing or approximating a Barabanov norm for a given set of matrices has consequently attracted some recent research interest \cite{Kozb1,Kozb2,MBara,Marg}. In this note we continue an investigation initiated in \cite{MBara} into the closely related question of determing when Barabanov norms are unique. Clearly if a given norm satisfies \eqref{barab}, then any positive scalar multiple of that norm will also satisfy \eqref{barab},  so when saying that a set of matrices has a ``unique'' Barbanov norm, we shall always mean only that any two Barabanov norms for that set must be proportional to one another by a scalar constant.

In the earlier article \cite{MBara} we established a sufficient condition for a finite irreducible set $\mathcal{A}=\{A_1,\ldots,A_m\}$ of $d \times d$ matrices over $\mathbb{R}$ or $\mathbb{C}$ to have a unique Barabanov norm in the sense defined above. We showed that if $\mathcal{A}$ has both the \emph{rank one property} and the \emph{unbounded agreements property}, defined formally in the next section, then a unique Barabanov norm for $\mathcal{A}$ exists. Roughly speaking, the unbounded agreements property states that there are not ``too many'' sequences $(A_{i_j})_{j=1}^\infty \in \mathcal{A}^{\mathbb{N}}$ such that the sequence of products $\varrho(\mathcal{A})^{-n}\|A_{i_n}\cdots A_{i_1}\|$ does not converge to zero in the limit as $n \to \infty$, while the rank one property states that for any fixed sequence of matrices, the vector space of all vectors $v$ such that $\varrho(\mathcal{A})^{-n}\|A_{i_n}\cdots A_{i_1}v\|$ converges to zero has the largest possible dimension. In this note we establish a new sufficient condition for the uniqueness of Barabanov norms which is complementary to the sufficient condition given in \cite{MBara}, and which applies in certain situations where there is instead a large supply of sequences $(A_{i_j})$ and vectors $v$ such that $\varrho(\mathcal{A})^{-n}\|A_{i_n}\cdots A_{i_1}v\|$ does not converge to zero. The new condition may also be applied to compact infinite sets of matrices. As direct examples of the application of the theorem, we exhibit firstly a pair of matrices which has a unique Barabanov norm but satisfies neither the rank one property nor the unbounded agreements property (hence not falling within the scope of \cite{MBara}), and secondly a countably infinite compact set of matrices which has a unique Barabanov norm, but such that every finite subset thereof has an uncountable family of Barabanov norms which are not proportional to one another.

We also use the main theorem in this note to investigate the robustness with respect to perturbation of the property of having a unique Barabanov norm. It was shown in \cite{MBara} that for every pair of integers $r,d \geq 2$, and for $\mathbb{K}$ equal to either $\mathbb{R}$ or $\mathbb{C}$, there exists an $r$-tuple of $d \times d$ matrices over $\mathbb{K}$ such that every sufficiently small perturbation of that $r$-tuple also has a unique Barbanov norm. 
As an application of our main theorem, we show that there exists a pair of real $2 \times 2$ matrices $\mathcal{A}$ with the following contrasting property: pairs of matrices having a unique Barbanov norm, and pairs of matrices not having a unique Barabanov norm, both form dense sets in a small open neighbourhood of $\mathcal{A}$. The property of having a unique Barabanov norm is thus shown to be highly sensitive to small perturbations of the set of matrices in certain cases.

\section{Statement and proof of main theorem}

Throughout the rest of this note we use the symbol $\mathbb{K}$ as a shorthand to denote either $\mathbb{R}$ or $\mathbb{C}$. Statements which are given in terms of $\mathbb{K}$ are thus valid if either of these two fields is consistently chosen. We use the symbol $\mathbf{M}_d(\mathbb{K})$ to denote the set of all $d \times d$ matrices over $\mathbb{K}$, which we equip with its usual topology as a normed vector space. The symbol $\|\cdot\|$ will be used to denote the Euclidean norm on $\mathbb{K}^d$, and also the corresponding induced matrix norm on $\Mat$. The symbol $\rho(B)$ will be used to denote the ordinary spectral radius of the matrix $B$.

If $\mathcal{A}$ is a compact subset of $\mathbf{M}_d(\mathbb{K})$ and $n \geq 1$ is an integer, we define $\mathcal{A}^n:=\{A_{i_1}\cdots A_{i_n} \colon A_i \in \mathcal{A}\}$. %The identity
%\[\varrho(\mathcal{A})=\inf_{n \geq 1} \sup\left\{\left\|A\right\|^{\frac{1}{n}} \colon A \in \mathcal{A}^n\right\} = \sup_{n \geq 1} \sup\left\{\rho\left(A\right)^{\frac{1}{n}} \colon A \in \mathcal{A}^n\right\}\]
%will occasionally be used to simplify calculations (check that we actually need it).
We shall say that $\mathcal{A}$ is \emph{product bounded} if the set $\bigcup_{n=1}^\infty \mathcal{A}^n$ is bounded. If $\varrho(\mathcal{A})$ is nonzero, we say that $\mathcal{A}$ is \emph{relatively product bounded} if $\varrho(\mathcal{A})^{-1}\mathcal{A}$ is product bounded. If $\mathcal{A}$ is irreducible then $\varrho(\mathcal{A})>0$ and $\mathcal{A}$ is relatively product bounded, see for example \cite{Jungers}. If $\mathcal{A} \subset \Mat$ is relatively product bounded, then following F. Wirth in \cite{Wirth} we define the \emph{limit semigroup} of $\mathcal{A}$ to be the set
\[\mathcal{S}(\mathcal{A}):=\bigcap_{m=1}^\infty\left(\overline{ \bigcup_{n=m}^\infty \varrho(\mathcal{A})^{-n}\mathcal{A}^n}\right) \subset \Mat.\]

We may now give the formal definition of the rank one property and the unbounded agreements property mentioned in the introduction. We say that $\mathcal{A}$ has the \emph{rank one property} if it is relatively product bounded and every nonzero element of $\mathcal{S}(\mathcal{A})$ is of rank one. The finite set of matrices $\mathcal{A}=\{1,\ldots,m\}$ has the \emph{unbounded agreements property} if for every pair of sequences $j_1,j_2 \colon \mathbb{N} \to \{1,\ldots,m\}$ such that $\limsup_{n \to \infty} \varrho(\mathcal{A})^{-n} \|A_{j_i(n)} \cdots A_{j_i(1)}\|>0$ for $i=1,2$, and every integer $\ell \geq 1$, there exist $k_1,k_2\geq 0$ such that $j_1(k_1+t)=j_2(k_2+t)$ for all $t$ in the range $1 \leq t \leq \ell$. (We do not define the unbounded agreements property for infinite sets of matrices.)

The central result of this note is the following sufficient condition for the uniqueness of the Barabanov norm:
\begin{theorem}\label{only}
Let $\mathcal{A}$ be a bounded, irreducible nonempty subset of $\Mat$ such that the limit semigroup $\mathcal{S}(\mathcal{A})$ has the following transitivity property: for every pair of nonzero vectors $v_1,v_2 \in \mathbb{K}^d$, there exist $B_1,B_2 \in \mathcal{S}(\mathcal{A})$ and $\lambda \in \mathbb{K}$ such that $B_1v_1=\lambda v_2$ and $B_2v_2 = \lambda^{-1}v_1$. Then $\mathcal{A}$ has a unique Barabanov norm.
\end{theorem}
\begin{proof}
Since $\mathcal{A}$ is irreducible, it admits at least one Barabanov norm.  Fix a nonzero vector $v_0 \in \mathbb{K}^d$ for the remainder of the proof, and suppose that $\vvv \cdot \vvv_1$ and $\vvv \cdot \vvv_2$ are both Barabanov norms for $\mathcal{A}$ which give norm $1$ to the vector $v_0$. To prove the theorem it is necessary and sufficient to show that $\vvv \cdot \vvv_1$ must be equal to $\vvv \cdot \vvv_2$.

Let $u \in \mathbb{K}^d$ be any nonzero vector. Choose $B_1, B_2 \in \mathcal{S}(\mathcal{A})$ and $\lambda \in \mathbb{K}$ such that $B_1v_0 = \lambda u$ and $B_2 u = \lambda^{-1}v_0$. Since $B_1,B_2 \in \mathcal{S}(\mathcal{A})$, there exist strictly increasing sequences of natural numbers $(n_{1,j})$, $(n_{2,j})$ and sequences of matrices $(A_{1,j})$, $(A_{2,j})$ such that $A_{i,j} \in \mathcal{A}^{n_{i,j}}$ for each $i=1,2$ and all $j \geq 1$, and $\lim_{j \to \infty} \varrho(\mathcal{A})^{-n_{i,j}}A_{i,j} = B_i$ for $i=1,2$. Since $\vvv \cdot \vvv_1$ is a Barabanov norm for $\mathcal{A}$, for each $j \geq 1$ we have
\[\vvv v_0 \vvv_1 = \max\left\{\varrho(\mathcal{A})^{-n_{1,j}} \vvv A v_0 \vvv_1 \colon A \in \mathcal{A}^{n_{1,j}}\right\} \geq \varrho(\mathcal{A})^{-n_{1,j}}\vvv A_{{1,j}}v_0\vvv_1 \]
and similarly
\[\vvv u \vvv_1 = \max\left\{\varrho(\mathcal{A})^{-n_{2,j}} \vvv A u \vvv_1 \colon A \in \mathcal{A}^{n_{2,j}}\right\} \geq \varrho(\mathcal{A})^{-n_{2,j}}\vvv A_{{2,j}}u\vvv_1.\]
Taking the limit as $j \to \infty$ in both of these inequalities and combining the results yields
\[1=\vvv v_0 \vvv_1 \geq  \vvv B_1 v_0 \vvv_1 = \vvv \lambda u \vvv_1 =|\lambda|.\vvv u \vvv_1 \geq |\lambda|.\vvv B_2u\vvv_1 = \vvv v_0\vvv_1=1\]
and therefore $\vvv u \vvv_1 = |\lambda|^{-1}$. However, we have not used any properties specific to $\vvv \cdot \vvv_1$ other than its being a Barabanov norm such that $\vvv v_0 \vvv_1=1$,  and so the above argument applies in an identical manner to the norm $\vvv \cdot \vvv_2$. We conclude that $\vvv u \vvv_1 = \vvv u \vvv_2$, and since the vector $u$ is arbitrary we have $\vvv \cdot \vvv_1 = \vvv \cdot \vvv_2$ as required. The proof is complete.
\end{proof}
\emph{Remark}. When $\mathbb{K}=\mathbb{R}$, the most straightforward case in which Theorem \ref{only} may be applied is that in which $\mathcal{S}(\mathcal{A})$ contains the special orthogonal group $SO(d)$, or more generally, when $\mathcal{S}(\mathcal{A})$ is simultaneously similar to a semigroup which contains $SO(d)$. In particular, if $\varrho(\mathcal{A})^{-1}\mathcal{A}$ contains a collection of matrices which generate a dense subsemigroup of $SO(d)$ (or which are simultaneously similar to a such a collection) then Theorem \ref{only} may be applied and $\mathcal{A}$ has a unique Barabanov norm. Similar remarks apply to the case $\mathbb{K}=\mathbb{C}$ and the group $SU(d)$. However, these cases certainly do not exhaust the possibilities of the theorem: for example, if $\mathcal{A}$ consists precisely of the set of rank one orthogonal projections on $\mathbb{R}^2$, then $\mathcal{S}(\mathcal{A})$ contains every real matrix which is equal to the composition of a rotation and an orthogonal projection, and Theorem \ref{only} also applies. Higher-dimensional examples of this type may of course also be constructed. In any case, Theorem \ref{only} is powerful enough to produce some interesting applications, which we describe in the following two sections.

\section{Examples}

In this section we provide some examples of the direct application of Theorem \ref{only}, before moving on to a theoretical application in the following section. In both examples we restrict our attention to the case $\mathbb{K}=\mathbb{R}$. The first example below illustrates the contrast between Theorem \ref{only} and \cite[Theorem 1.5]{MBara}.

{\bf{Example 1.}} Define
\[A_1:=\left(\begin{array}{cc}1&0\\0&0\end{array}\right), \qquad A_2:=\left(\begin{array}{cc}\cos \theta\pi & -\sin \theta\pi\\\sin \theta\pi & \cos\theta\pi\end{array}\right)\]
and $\mathcal{A}:=\{A_1,A_2\}$, where $\theta \in \mathbb{R} \setminus \mathbb{Z}$. Then $\mathcal{A}$ is irreducible with $\varrho(\mathcal{A})=1$, has neither the rank one property nor the unbounded agreements property, and has a unique Barbanov norm if and only if $\theta \notin \mathbb{Q}$. 
\begin{proof}
Let us first establish the properties of $\mathcal{A}$ for general $\theta \in \mathbb{R} \setminus \mathbb{Z}$. Since $\theta \notin \mathbb{Z}$ the matrix $A_2$ does not preserve any one-dimensional subspace of $\mathbb{R}^2$, and therefore $\mathcal{A}$ is irreducible. It is straightforward to see that $\max\{\|A\| \colon A \in \mathcal{A}^n\}=1$ for every $n \geq 1$ and consequently $\varrho(\mathcal{A})=1$. In particular $\mathcal{A}$ is product bounded. Every accumulation point at infinity of the sequence $( A_2^n)_{n=1}^\infty$ has rank two, and so $\mathcal{A}$ does not have the rank one property. Since $\lim_{n \to \infty}\|A_i^n\|=1$ for both $i=1$ and $i=2$ the unbounded agreements property is also not satisfied.

Let us now consider those properties which depend on whether or not $\theta \in \mathbb{Q}$. If $\theta \notin \mathbb{Q}$ then every rotation matrix in $\mathbf{M}_2(\mathbb{R})$ is an accumulation point of $(A_2^n)_{ n = 1}^\infty$, and so $\mathcal{S}(\mathcal{A})$ contains the group of rotation matrices. It follows easily that $\mathcal{A}$ meets the hypotheses of Theorem \ref{only} and has a unique Barabanov norm. Conversely, let us suppose that $\theta \in \mathbb{Q}$. Let $\vvv\cdot\vvv$ be any norm whose unit ball is invariant under rotation through angle $\theta\pi$ and such that $\vvv (x,0)^T \vvv \leq \vvv(x,y)^T\vvv$ for all $x,y \in \mathbb{R}$. The former property ensures that $\vvv A_2v\vvv = \vvv v \vvv = \varrho(\mathcal{A})\vvv v \vvv$ for every $v \in \mathbb{R}^2$, and the latter property ensures that $\vvv A_1v\vvv \leq \vvv v \vvv = \varrho(\mathcal{A})\vvv v \vvv$ for every $v$, so in particular any such norm is Barabanov. If $K \subset \mathbb{R}^2$ is a compact convex set with nonempty interior which is symmetrical with respect to rotation about the origin through angles $\theta\pi$ and $\pi$, and such that there exists a vertical tangent to $K$ at all of its boundary points which lie on the horizontal axis, then $K$ is the unit ball of a norm which has the required properties. It is clear that uncountably many such sets exist which are not related to one another by scalar multiplication, and we conclude that $\mathcal{A}$ has uncountably many Barabanov norms. 
\end{proof}

The following example shows that the uniqueness of Barbanov norms can be a quite delicate phenomenon:

{\bf{Example 2.}} Let us define a compact subset of $\mathbf{M}_2(\mathbb{R})$ by
\[\mathcal{A}:=\left\{\left(\begin{array}{cc}1&0\\0&1\end{array}\right)\right\}\cup \left\{\left(\begin{array}{cc}\cos \frac{\pi}{2^n} & -\sin  \frac{\pi}{2^n}\\ \sin \frac{\pi}{2^n} & \cos \frac{\pi}{2^n}\end{array}\right)\colon n \geq 1\right\}.\]
Then every nonempty subset of $\mathcal{A}$ has joint spectral radius equal to one, is product bounded, and does not have the rank one property. Every subset of $\mathcal{A}$ with at least two elements is irreducible, and finite subsets of $\mathcal{A}$ which have at least two elements do not satisfy the unbounded agreements propery. Every infinite subset of $\mathcal{A}$ has a unique Barabanov norm, but every finite nonempty subset of $\mathcal{A}$ has uncountably many Barabanov norms.
\begin{proof}
Let $\mathcal{B} \subseteq \mathcal{A}$ be a nonempty subset. It is clear that $\sup\{\|B\|\colon B \in \mathcal{B}^n\}=1$ for all $n \geq 1$ so that $\varrho(\mathcal{B})=1$ and $\mathcal{B}$ is product bounded, and it is also clear that $\mathcal{S}(\mathcal{B})$ contains the identity matrix so that the rank one property is not satisfied. If $\mathcal{B}$ has at least two elements then it includes a rotation matrix with no real eigenvalues, and hence $\mathcal{B}$ is irreducible. Since $\lim_{n \to \infty}\|B^n\| =1$ for every $B \in \mathcal{B}$, the unbounded agreements property is not satisfied when $\mathcal{B}$ is finite and contains at least two elements. Since every element of $\mathcal{B}$ preserves the Euclidean norm, that norm is a Barabanov norm for $\mathcal{B}$.

Let us now consider the uniqueness or otherwise of Barabanov norms for $\mathcal{B}$. Suppose first that $\mathcal{B}$ is infinite. In this case there exist infinitely many positive integers $q$ such that $\mathcal{B}$ includes the matrix of rotation through angle $\pi/2^q$. If $q$ is such an integer, then in particular it follows that $\mathcal{S}(\mathcal{B})$ contains the group of all rotations through angles of the form $k\pi/2^q$. Since $q$ may be taken arbitrarily large it follows that $\mathcal{S}(\mathcal{B})$ contains every rotation by a dyadic rational multiple of $\pi$, and since $\mathcal{S}(\mathcal{B})$ is closed we conclude that $SO(2)\subseteq \mathcal{S}(\mathcal{B})$. Theorem \ref{only} therefore applies and the Euclidean norm is the unique Barbanov norm of $\mathcal{B}$.

Now let us suppose that $\mathcal{B} \subset  \mathcal{A}$ is finite and nonempty. If $\mathcal{B}$ consists only of the identity matrix then every norm on $\mathbb{R}^2$ is preserved by $\mathcal{B}$ and hence is Barabanov. Otherwise, let $n$ be the largest integer such that $\mathcal{B}$ contains the matrix corresponding to rotation through angle $\pi/2^n$. If $K \subset \mathbb{R}^2$ is a compact convex set with nonempty interior which is invariant under rotation through angle $\pi/2^n$, then it is invariant under the action of every element of $\mathcal{B}$. To each such $K$ there corresponds a norm on $\mathbb{R}^2$  which has $K$ as its unit ball, hence is invariant under every element of $\mathcal{B}$ and therefore is Barabanov. Since there exist uncountably many compact convex sets $K$ which are invariant under rotation through angle $\pi/2^n$ and are not pairwise similar, we conclude that $\mathcal{B}$ has uncountably many Barabanov norms. \end{proof}

\section{A theoretical application}

Following the notation of \cite{MBara}, let us use the symbol $\mathcal{O}_2(\mathbb{R}^2)$ to denote the set of all ordered pairs of $2 \times 2$ real matrices, which we equip with the topology arising from the natural identification of this space with $\mathbf{M}_2(\mathbb{R}) \oplus \mathbf{M}_2(\mathbb{R})$. In \cite{MBara} we showed that $\mathcal{O}_2(\mathbb{R}^2)$ contains a nonempty open set $\mathcal{U}$ with the property that for every $(A_1,A_2) \in \mathcal{U}$, the set $\mathcal{A}=\{A_1,A_2\}$ has a unique Barabanov norm. For matrix pairs belonging to $\mathcal{U}$, therefore, the property of having a unique Barabanov norm is robust with respect to sufficiently small perturbations of either or both of the matrices comprising the pair. This result naturally leads one to ask whether this phenomenon is typical: is the set of all $(A_1,A_2) \in \mathcal{O}_2(\mathbb{R}^2)$ such that $\{A_1,A_2\}$ has a unique Barabanov norm open, or does it contain an open set which is dense in $\mathcal{O}_2(\mathbb{R}^2)$? The following result answers both questions negatively, and shows that the uniqueness of Barabanov norms can be sensitive to small perturbations in quite a strong manner:
\begin{theorem}
There exists a nonempty open set $\mathcal{V} \subset \mathcal{O}_2(\mathbb{R}^2)$ such that the sets
\[\mathcal{V}_1:=\left\{(A_1,A_2) \in\mathcal{V}\colon \{A_1,A_2\}\text{ has a unique Barabanov norm}\right\}\]
and
\[\mathcal{V}_2:=\left\{(A_1,A_2) \in\mathcal{V}\colon \{A_1,A_2\}\text{ does not have a unique Barabanov norm}\right\}\]
are both dense in $\mathcal{V}$.
\end{theorem}
\begin{proof}
Let $A_2$ be a rotation matrix which does not have real eigenvalues, and let $A_1$ be any matrix such that $\|A_1\|<\|A_2\|=1$. We will take $\mathcal{V}$ to be a suitably small neighbourhood of $(A_1,A_2)$.

For each $\delta>0$ let $\mathsf{B}_\delta$ denote the open ball about the origin in $\mathbf{M}_2(\mathbb{R})$ which has radius $\delta$ with respect to the spectral norm. Since the eigenvalues of $A_2$ are simple, we may choose $\varepsilon>0$ small enough that there exist continuous functions $E \colon \mathsf{B}_\varepsilon \to \mathbb{C}$, $V \colon \mathsf{B}_\varepsilon \to \mathbb{C}^2$ such that for all $C \in \mathsf{B}_\varepsilon$, $E(C)$ is a strictly complex eigenvalue of the real matrix $A_2+C$ with corresponding complex eigenvector $V(C)$. Since $A_2+C$ is real, the complex conjugates of $E(C)$ and $V(C)$ are also an eigenvalue and an eigenvector respectively. Since $E(C)$ is strictly complex it is not equal to $\overline{E(C)}$, and consequently the associated eigenvectors $V(C)$ and $\overline{V(C)}$ are linearly independent over $\mathbb{C}$. It follows from this that the real and imaginary parts of $V(C)$ are a linearly independent pair of vectors with respect to $\mathbb{R}$. For each $C \in \mathsf{B}_\varepsilon$ let us now define $S(C)$ to be the invertible real matrix with first column given by $\Im( V(C))$ and second column given by $\Re (V(C))$. An elementary calculation shows that
$S(C)^{-1}(A_2+C)S(C)$ is precisely the real matrix of rotation through angle $\arg E(C)$ multiplied by the positive scalar factor $|E(C)|=\rho(A_2+C)$. For each $C \in \mathsf{B}_\varepsilon$ define a norm $\|\cdot\|_C$ on $\mathbb{R}^2$ by $\|v\|_C:=\|S(C)^{-1}v\|$ for every $v \in \mathbb{R}^2$. It is easily seen that $\rho(A_2+C)^{-1}(A_2+C)$ is an isometry of $\mathbb{R}^2$ with respect to this norm, and in particular $\|A_2+C\|_C = \rho(A_2+C)$. Since $A_2$ is a rotation matrix, a direct calculation shows that $S(0)$ is proportional to the identity and therefore $\|\cdot\|_0$ is a scalar multiple of the Euclidean norm on $\mathbb{R}^2$.

Let us now define
\[\mathcal{V}:=\left\{(B_1,B_2) \in \mathcal{O}_2(\mathbb{R}^2) \colon B_2-A_2 \in \mathsf{B}_\varepsilon\text{ and }0<\|B_1\|_{B_2-A_2}<\rho(B_2)\right\}.\]
Clearly $\mathcal{V}$ contains $(A_1,A_2)$, and since $S \colon \mathsf{B}_\varepsilon \to \mathbf{M}_2(\mathbb{R})$ is continuous, $\mathcal{V}$ is open. We claim that $(B_1,B_2) \in \mathcal{V}$ has a unique Barabanov norm if and only if the eigenvalues of $\rho(B_2)^{-1}B_2$ are not roots of unity. An easy perturbation argument shows that pairs $(B_1,B_2)\in\mathcal{V}$ such that the eigenvalues of $B_2$ have irrational arguments, and pairs such that the eigenvalues of $B_2$ have rational arguments, are both dense in $\mathcal{V}$. It follows that establishing this claim is sufficient to complete the proof of the theorem. 

For the rest of the proof let us fix an arbitrary pair of matrices $(B_1,B_2) \in \mathcal{V}$. Define $C:=B_2-A_2 \in \mathsf{B}_\varepsilon$. It is straightforward to see that $\sup \{\|B\|_C \colon B \in \mathcal{B}^n\}=\rho(B_2)^n$ for every $n \geq 1$ and therefore $\varrho(\mathcal{B})=\rho(B_2)$. Since $\rho(B_2)^{-1}B_2$ is an isometry of $\mathbb{R}^2$ with respect to the norm $\|\cdot\|_C$, and $\|B_1\|_C<\rho(B_2)=\varrho(\mathcal{B})$, it follows directly that $\|\cdot\|_C$ is a Barabanov norm for $\mathcal{B}$.

Let us suppose first that the eigenvalues of $\rho(B_2)^{-1}B_2$ are not roots of unity. In this case $\rho(B_2)^{-1}S(C)^{-1}B_2S(C)$ is a matrix corresponding to rotation through an irrational angle. It follows that every rotation matrix is a limit point at infinity of the sequence $(\varrho(\mathcal{B})^{-n}S(C)^{-1}B_2^nS(C))_{n=1}^\infty$, and so $S(C)^{-1}\mathcal{S}(\mathcal{B})S(C)$ contains the group of rotation matrices. We deduce that $\mathcal{B}$ satisfies the hypotheses of Theorem \ref{only} and conclude that $\|\cdot\|_C$ is the unique Barabanov norm of $\mathcal{B}$.

Now let us suppose instead that the eigenvalues of $\rho(B_2)^{-1}B_2$ are roots of unity. Similarly to examples 1 and 2 in the previous section, there exist uncountably many norms on $\mathbb{R}^2$ which are preserved by the rational-angle rotation matrix $\rho(B_2)^{-1}S(C)^{-1}B_2 S(C)$ and are not proportional to one another or to the Euclidean norm. Modifying these norms by composition with $S(C)^{-1}$ as in the definition of the norm $\|\cdot\|_C$, we obtain an uncountable family of norms on $\mathbb{R}^2$ which are preserved by $\rho(B_2)^{-1}B_2$ and are not proportional to one another or to $\|\cdot\|_C$.  Let $\vvv \cdot \vvv$ be any such norm. We will show that for every sufficiently small real number $\kappa>0$, the norm on $\mathbb{R}^2$ given by $\vvv v \vvv_*:=\|v\|_C+\kappa \vvv v \vvv$ is a Barabanov norm for $\mathcal{B}$. By repeating this procedure using a different norm $\vvv \cdot \vvv$ which is also preserved by $\rho(B_2)^{-1}B_2$, or indeed by simply varying the constant $\kappa$ within its permitted range, it is clear that we may obtain an uncountable family of Barabanov norms for $\mathcal{B}$ which are not pairwise proportional to one another. 

Given a norm $\vvv \cdot \vvv$ on $\mathbb{R}^2$ which is preserved by $\rho(B_2)^{-1}B_2$, then, let us define $\xi:=\rho(B_2)^{-1}\|B_1\|_C \in (0,1)$, choose a constant $K>1$ such that $K^{-1}\|v\|_C \leq \vvv v \vvv \leq K\|v\|_C$ for all $v \in \mathbb{R}^2$, and let $\kappa>0$ be small enough that $1+K \kappa<\xi^{-1}$. Define $\vvv v\vvv_*:=\|v\|_C+\kappa\vvv v \vvv$ for all $v \in \mathbb{R}^2$. Clearly $\vvv B_2v\vvv_* = \rho(B_2)\vvv v \vvv_*= \varrho(\mathcal{B})\vvv v \vvv_*$ for every $v \in \mathbb{R}^2$, and so to prove that $\vvv \cdot \vvv_*$ is a Barabanov norm for $\mathcal{B}$ it is sufficient to show that $\vvv B_1\vvv_* \leq \rho(B_2)=\varrho(\mathcal{B})$. For each $v \in \mathbb{R}^2$ we have
\begin{align*}
\vvv B_1v  \vvv_* &= \|B_1v\|_C +\kappa \vvv B_1v\vvv\\
 &\leq (1+K\kappa)\|B_1v\|_C\\
& \leq \xi\rho(B_2)(1+K\kappa)\|v\|_C\\
&=\xi\rho(B_2)\left(1+\left(K-\frac{1}{K}\right)\kappa\right)\|v\|_C  +\frac{\xi \rho(B_2)\kappa}{K}\|v\|_C\\
&< \rho(B_2)\|v\|_C + \xi\rho(B_2)\kappa\vvv v \vvv\\
& <\rho(B_1)\vvv v \vvv_*,
\end{align*}
and it follows that $\vvv B_1 \vvv_*<\rho(B_2)$ so that $\vvv \cdot \vvv_*$ is a Barabanov norm for $\mathcal{B}$ as claimed.  This proves the claim and hence completes the proof of the theorem.
\end{proof}

\section{Acknowledgment}

This research was conducted as part of the ERC grant MALADY (246953).

\bibliographystyle{siam}
\bibliography{2Barabanov}

\end{document}